\documentclass[12pt]{scrartcl}

 \usepackage{amsmath,amsxtra,amscd,amssymb,latexsym,stmaryrd,amsthm}
 \usepackage{accents, mathrsfs}
  \usepackage[utf8]{inputenc}
  \usepackage[T1]{fontenc}

 \usepackage{authblk}

\addtokomafont{title}{\rmfamily\itshape}

\theoremstyle{plain}
\newtheorem{theo}{Theorem}[section]
\newtheorem*{theo*}{Theorem}

\theoremstyle{definition}

\newtheorem{rema}[theo]{Remark}

\newtheorem{ques}[theo]{Question}

\newcommand*{\dd}%
  {\relax\ifnum\lastnodetype>0\mskip\medmuskip\fi\mathrm{d}}

\newcommand{\Lie}{\mathscr{L}}

\newcommand{\derivative}[1]{\left.\frac{\dd}{\dd #1}\right|_{#1=0}}
\newcommand{\bdot}[1]{\accentset{\mbox{\large\bfseries .}}{#1}}

\title{The linear request problem}

 \author[$\dag$]{Benoit R. Kloeckner\thanks{Supported by the Agence Nationale de la Recherche, grant ANR-11-JS01-0011.}}

 \affil[$\dag$]{{\small Universit\'e Paris-Est, Laboratoire d'Analyse et de Mat\'ematiques Appliqu\'ees (UMR 8050), UPEM, UPEC, CNRS, F-94010, Cr\'eteil, France}}

\begin{document}

\maketitle

\begin{abstract}
We propose a simple approach to a problem
introduced by Galatolo and Pollicott, consisting in perturbing a
dynamical system in order for its absolutely continuous invariant
measure to change in a prescribed way. Instead of using transfer operators, we observe that restricting to an infinitesimal conjugacy already yields a solution. This allows us to work in any dimension and dispense from any dynamical hypothesis. In particular, we don't need to assume hyperbolicity to obtain a solution, although expansion moreover ensures the existence of an infinite-dimensional space of solutions.
\end{abstract}

Let $M$ be a compact $n$-manifold and $T:M\to M$ be a smooth map, seen as a discrete-time dynamical system. Assume that there is a positive smooth invariant measure $\omega$, also seen as a nowhere vanishing $n$-form.

The \emph{linear response theory}(see e.g. \cite{Ruelle2,BS2}) is typically concerned with the following question: if $(T_t)_{t\in(-\varepsilon,\varepsilon)}$ is a family of maps with $T_0=T$, differentiable at $t=0$, such that each $T_t$ preserves a smooth measure $\omega_t$, can we differentiate $\omega_t$ with respect to $t$ and express $\bdot{\omega}_0 := \derivative{t}\omega_t$ in term of $\bdot{T}_0 :=\derivative{t} T_t$?

A recent article by Galatolo and Pollicott \cite{GP} proposes to study the opposite direction; they roughly ask two questions (and provide some answers for expanding maps of the circle), the second of which we voluntarily keep imprecise:
\begin{ques}\label{ques:1}
Given a smooth function $\rho:M\to\mathbb{R}$ of vanishing integral,\footnote{The condition $\int_M\rho \omega=0$ is needed whenever invariant measures are normalized to have fixed total mass.} can we find a perturbation $(T_t)_{t\in(-\varepsilon,\varepsilon)}$ of $T$, differentiable at $t=0$ and preserving a family of smooth measure $(\omega_t)_{t\in(-\varepsilon,\varepsilon)}$ such that $\bdot{\omega}_0 = \rho \omega$? Can we then express the possible values of $\bdot{T}_0$ in term of $\rho$?
\end{ques}

\begin{ques}\label{ques:2}
When the previous question has a positive answer, can we find an ``optimal'' value of $\bdot{T}_0$?
\end{ques}

The goal of this note is to observe that Question \ref{ques:1} has a positive answer in a very general setting, without dynamical hypothesis and in every dimension. This follows from well-known facts in differential geometry, and the only faint bit of novelty is in the observation that one can restrict to conjugate deformations, of the form $T_t=\varphi^{t} \circ T \circ \varphi^{-t}$ where $(\varphi^t)_{t\in(-\varepsilon,\varepsilon)}$ is the flow of a vector field.\\


Let us first define properly what it means for a family $(T_t)_{t\in(-\varepsilon,\varepsilon)}$ to be differentiable at $t=0$. Galatolo and Pollicott work on the circle and implicitly use its parallelism (all tangent spaces can be identified), which is not possible on a general manifold. Pointwise, we want to ask that for each $x\in M$, the curve $(T_t(x))_{t\in(-\varepsilon,\varepsilon)}$ is differentiable at $t=0$, and we want the derivative of this curve to depend smoothly on $x$. We want to stress that $\bdot{T_0}(x):=\derivative{t} T_t(x)$ is an element of $T_{T(x)} M$, not of $T_xM$, and $\bdot{T_0}$ is thus not a vector field.
We will thus consider the set $\Gamma_T(M)$ of smooth maps $Z : M \to TM$ such that for all $x\in M$, $Z_x\in T_{T(x)} M$, and say that
a family $(T_t)_{t\in(-\varepsilon,\varepsilon)}$ of smooth maps $M\to M$ is differentiable at $t=0$ if $\bdot{T_0}(x)$ is defined for all $x$ and $\bdot{T_0}\in\Gamma_T(M)$.

For example, if $X$ is a smooth vector field on $M$, we can consider its flow $(\varphi^t)_{t\in\mathbb{R}}$ and the family $T_t := \varphi^{t} \circ T \circ \varphi^{-t}$. Then a direct computation shows that $(T_t)_{t\in\mathbb{R}}$ is differentiable at $0$:
\begin{align*}
\varphi^t\circ T\circ \varphi^{-t}(x)
  &= \varphi^t\circ T\big(x-tX_x+o(t)\big) \\
  &= \varphi^t\Big(T(x)-t D_xT(X_x)+o(t)\Big) \\
  &= T(x)-t D_xT(X_x)+tX_{T(x)}+o(t)
\end{align*}
so that $\bdot{T_0}(x)=-D_xT(X_x)+X_{T(x)}$, which we also write $\bdot{T_0}=-DT(X)+X_T$. This is naturally an element of $\Gamma_T(M)$, as it should.

Observe that since $T$ preserves $\omega$, $T_t$ then preserves $\omega_t := \varphi^{t}_*\omega$. But
\[\derivative{t} \varphi^{t}_*\omega = -\Lie_X\omega\]
where $\Lie$ denotes the Lie derivative.
To answer Question \ref{ques:1} by restricting to such deformations by conjugacy, we thus only have to check that for all $\rho$ such $\int_M \rho \omega =0$, there is a vector field $X$ such that $\Lie_X \omega = -\rho\omega$. This is well-known, see e.g. the proof of Moser's theorem in \cite{KH} (Theorem 5.1.27 page 195), but let us recall the argument for the sake of completeness.

Since $\omega$ is a $n$-form, the Cartan formula reads
$\Lie_X\omega = \dd(i_X\omega)$ where $i_X\omega$  is the $(n-1)$-form $\omega(X; \cdot;\cdots;\cdot)$ obtained by contraction with $X$. 
Since $M$ is compact, its top-dimensional cohomology is $1$-dimensional, generated by the class of $\omega$ (or any volume form), and every $n$-form of vanishing total integral is exact. This means that there exist a $(n-1)$ form $\theta$ such that $d\theta=-\rho\omega$. Now, $\omega$ being a volume form it is non-degenerate, which means that any $(n-1)$-form can be obtained by contracting $\omega$ with a well-chosen vector field; in particular, there must exist a vector field $X$ such that $i_X\omega=\theta$. Putting all this together, we get the desired conclusion:
\[\bdot{\omega}_t = -\Lie_X\omega = -\dd(i_X\omega)=-d\theta=\rho\omega.\]

We conclude:
\begin{theo}\label{theo:main}
Let $T:M\to M$ be a smooth map acting on a compact smooth Riemannian manifold, preserving a smooth volume form
$\omega$. Let $\rho:M\to\mathbb{R}$ be a smooth function
such that $\int_M \rho \omega = 0$.

There exist a deformation $(T_t)_{(-\varepsilon,\varepsilon)}$ of $T$ that is differentiable at $t=0$ and preserving smooth volume forms $\omega_t$, such that $\bdot{\omega}_0=\rho\omega$. Moreover one can ask all $T_t$ to be smoothly conjugate to $T$.
\end{theo}

Now, observe that in the above $\theta$ is not unique: one can add to it any closed $(n-1)$ form, and will obtain a different suitable $X$. If $n\ge2$, the space of closed $(n-1)$ form is infinite-dimensional (it contains all exact forms $\dd\alpha$ where $\alpha$ is a $(n-2)$-form); but if $M=\mathbb{S}^1$, there are no $(n-2)$-forms and the space of closed $(n-1)$-form is the space of functions with vanishing derivative, i.e. of constant functions.

This lack of uniqueness points, in dimension at least $2$, to the existence of many solutions to Question \ref{ques:1} even when restricting to conjugate deformations; however, since we consider perturbations only at their first order, we have to determine whether the space of possible $\bdot{T_0}$ is indeed as large as the set of vector fields $X$ generating the conjugacy.

We thus consider smooth vector fields $X$ and $Y$  on $M$, such that $\Lie_X\omega = \Lie_Y\omega=-\rho\omega$ and
$-DT(X)+X_T = -DT(Y)+Y_T$. Does this imply that $X=Y$? The answer strongly depends on $T$.
The first equation means that $X-Y$ preserves $\omega$, while the second equation can be rewritten as
\[X_T-Y_T = DT(X-Y),\]
i.e. the vector field $X-Y$ is $T$-invariant. The difference between the spaces of possible values for $X$ and for $\bdot{T_0}$ above is therefore determined by the space of $T$-invariant vector fields preserving $\omega$. There is no general rule: there could be no non-zero such vector fields (e.g. if $T$ is expanding) or there could be plenty (e.g. if $T$ is the identity, which of course is a very degenerate case: we get $\bdot{T_0}=0$). Depending on this, there could be an infinite-dimensional space of possible values for $\bdot{T_0}$ in Theorem \ref{theo:main}, or a unique one, or anything in between.

We conclude with a few remarks.
\begin{rema}
Sometimes one consider perturbations of the form $\bdot{T_0}(x)=X_{T(x)}$ for some vector field $X$. This is equivalent to considering $\bdot{T_0}\in\Gamma_T(M)$ if $T$ is invertible, but is less general otherwise as the images of any $x,y\in M$ such that $T(x)=T(y)$ would be asked to be perturbed identically. \end{rema}

\begin{rema}
Using Moser's theorem \cite{Moser} instead of its first-order version, one sees that for all volume form $\omega'$, there is a conjugate of $T$ that preserves $\omega'$.
\end{rema}

\begin{rema}
It could be asked what happens in lower regularity, e.g. $C^k$. We do not enter into these details, since the strategy would be the very same, only keeping track of the available regularity for the various objects $\omega$, $\theta$, $X$, etc.
\end{rema}

\begin{rema}\label{rema:q2}
One can wonder whether the observation here can shed some light on Question \ref{ques:2}. If one measures optimality by fixing a Riemannian metric on $M$ and trying to minimize $\lVert \bdot{T_0} \rVert_{L^2}$, as Galatolo an Pollicott do in \cite{GP}, then in general the optimal first-order perturbation is not of the form obtained here by conjugacy. This is easily seen in the case of the doubling map on the circle.

One can also restrict to conjugate perturbations, and try to minimize $\lVert X\rVert_{L^2}$ where $X$ is the vector field generating the conjugating diffeomorphisms. Then the problem directly relates to solving the modified Poisson equation
\[\nabla\cdot (\eta\nabla u) = -\rho\]
where $\eta,\rho$ are the densities of the invariant measure and of its perturbation with respect to the Riemannian volume and $u$ is the unknown (then $\nabla u$ is the optimal vector field $X$). This is an elliptic PDE, whose theory is well established and that admits smooth solutions. The first version of this note used this to obtain Theorem \ref{theo:main}, but not insisting on finding a \emph{gradient} vector field can be done even more easily, as above. I wish to thank Anatole Katok for interesting criticism on the first version of this note, pointing to this simplification.
\end{rema}

\bibliographystyle{alpha}
\bibliography{LR}
\end{document}